\documentstyle[11pt]{article}

\newtheorem{theorem}{{Theorem}}[section]
\newtheorem{proposition}[theorem]{{Proposition}}%[section]
\newtheorem{definition}[theorem]{{Definition}}%[section]
\newtheorem{corollary}[theorem]{{Corollary}}%[section]
\newtheorem{fact}[theorem]{{Fact}}%[section]
\newtheorem{remark}[theorem]{{\it Remark}}%[section]
\newtheorem{step}{{Step}}

\newenvironment{proof}{{\it Proof.}}{\hfill$\Box$\medskip}

\begin{document}

%macros

%titre

\title{Isometry groups and geodesic
foliations of Lorentz manifolds. \hbox{Part II}:
Geometry of analytic Lorentz manifolds with large isometry groups}
\author{Abdelghani Zeghib}
 %\date{November 24, 1996}
\date{\today}
\maketitle

\begin{abstract} This is Part  II of
a series on  noncompact isometry groups of Lorentz
manifolds. We have introduced in  Part I, a compactification
of these isometry groups, and called ``bi-polarized'' those
Lorentz manifolds having a ``trivial '' compactification.
Here we show a geometric rigidity of non-bi-polarized
Lorentz manifolds; that is,  they are (at least locally)
warped products of constant curvature Lorentz manifolds by
Riemannian manifolds.

\end{abstract}

%\tableofcontents

\section{Introduction}

We continue here our investigation of noncompact
isometry groups of compact Lorentz manifolds, started
in Part I \cite{Part.I}
%Zeg.Lorentz.dynamics}.
%The part I
%presented
which contains  dynamical ingredients.
Its fundamental tool was the notion of
approximate stability. This second  part (which is
in fact fairly independent of  Part I) is
geometrical, and has the {\bf warped
product} construction as
a fundamental tool. Recall
that this is a construction
in the class of
pseudo-Riemannian manifolds,
defined as follows.
Let $(L, h)$ and
$(N, g)$ be  two pseudo-Riemannian
  manifolds, and $ w: L \to {\bf R}^+$
a (warping) function. The warped
product $M= L\times_wN$, is the topological
product
  $L \times N$, endowed with
the  metric
$h \bigoplus w g$.

 The  warped product
 construction is very useful
in Riemannian  as well as Lorentzian geometry, since it gives
%offers
%  permit to get
%elaborated
sophisticated examples  from simple ones.
 For instance,  warped product models  are omnipresent in
cosmological theories.

Here we are interested in the case where  $L$ is Riemannian,
$N$ is Lorentzian, and hence $M = L \times_wN$
is Lorentzian.
However, from a physical viewpoint, it is
more
interesting
%realistic
to consider the situation where
$L$ is Lorentzian, and
$N$ is Riemannian (the warping function is
thus a universe expansion function).

There are two key properties of warped products. \\

1) If $f: N \to N$ is an isometry, then the
trivial extension $\bar{f}: (x, y) \in
L \times N \to (x, f(y)) \in L\times N$, is
an isometry of
$L\times_w N$ (see \S \ref{extension}).

In particular, in the class of Lorentz manifolds with
large isometry groups, one can perform
warped products by (any) Riemannian manifolds.

In fact,  the warped products  are reminiscent of
semi-direct products in the category of groups,
the factor $N$ playing the role ofthe normal subgroup.
One
may  justify this by the fact that, indeed,
Isom$(N)$ is a normal subgroup of the subgroup
of elements of Isom$(L \times_w N)$, which preserve
the topological product $L \times N$
(i.e. the foliation determined by the factors
$L$ and $N$). This suggests to us to
call
the factor $N$ the {\bf normal} factor
of the warped product.
 (This  will be useful for us because
we actually need to to
distinguish between the factors).  \\

2) The second fundamental
fact about warped products is that if $S$ is a geodesic submanifold of
$N$, then $L \times S$
 is a geodesic submanifold in $L \times_wN$
(see \ref{geodesic.submanifold}).

It is very special when  a
Lorentz manifold (or in general a pseudo-Riemannian manifold,
or even just a manifold endowed
with a connection) admits many geodesic submanifolds of
dimension $>1$ (and codimension $\neq 0$).
Generically, there is no such submanifold.
The degenerate case, when every tangent plane
is tangent to a geodesic submanifold, corresponds exactly
to Lorentz manifold of constant curvature  (this is also  true
in the general pseudo-Riemannian case), see \S \ref{tautologic}.

Here, we are especially interested in
the case where the factor $N$ has  constant curvature. So, like
$N$, $L \times_wN$ has many geodesic hypersurfaces.  \\

In this article we investigate the
relationships
between the following three phenomena: being a warped
 product, having a large isometry
group, and having abundant geodesic hypersurfaces.
In particular, we show that in some situations, one of the
%sub-properties
%(1) or (2)
second two properties
may lead to a warped product structure.

\subsection{Abundance of geodesic hypersurfaces leads
to a warped product structure}
\label{introduction.def.1}

Let $M$ be  Lorentz manifold. In order to
analyze
the set of geodesic hypersurfaces in
$M$, we associate to
any $x \in M$ a set $C_x$ of tangent directions
at $x$  (i.e. 1-dimensional sub spaces of $T_xM)$ defined
as follows.
 A direction
$u \in {\bf P}(T_xM)$ belongs to
 $C_x$, if it is {\bf isotropic}
(this choice is related to our anti-physical preference
of signatures of the factors $L$ and
$N$), and
 the orthogonal $u^\perp$
 determines a geodesic hypersurface.
 Equivalently,  there
is a {\bf lightlike} geodesic hypersurface
$H$ passing through $x$, such that $T_xH = u^\perp$.
(Recall  here that, if
$<,>$ denotes the Lorentz scalar on $T_xM$,
then a vector $v$ is isotropic, if
$<v,v> =0$, and a hyperplane $E \subset
T_M$ is lightlike,
if its orthogonal is isotropic, or equivalently, the
restriction of $<,>$ on $E$ is degenerate. A hypersurface
is lightlike if its tangent space is everywhere lightlike).

Consider  the open set ${\cal W}(M)$ of  points $M$
having a neighborhood isometric to a warped product
with the normal factor being
 a Lorentz manifold of constant curvature and dimension
$\geq 3$.
%, by some Riemannian manifold.
That is, $x \in
{\cal W}(M)$, if and only if there is a neighborhood $U$ of
$x$ isometric to a warped
product
$L \times_w N$, where
$N$ is a Lorentz manifold of  constant curvature and
dim $N \geq 3$.

Section \ref{tautologic} is devoted
to the study of the relationship
 between the map,
$x \to C_x$, and
${\cal W}(M)$. The following statement is
a simple corollary of this study, which will be
fully proved only in the analytic case, but it
needs some results from \cite{Zeg.tautologic}
 in the smooth case.
%(see \cite{Zeg.tautologic}).

\begin{theorem} Let $M$ be a Lorentz manifold. Suppose
that $C_x$ is infinite for
all
$x \in M$. Then,
${\cal W}(M)$ is dense (and open by definition) in
$M$.
\end{theorem}

Observe that this a local result and that
its  converse is obviously true.
This result admits a kind of generalization
to the general pseudo-Riemannian case. Notice
that the condition on the existence of
geodesic hypersurfaces, cannot be relaxed to an existence
condition of geodesic submanifolds
of higher codimension.
For instance, in the Riemannian case, the symmetric space
${\bf C}P^n$ admits many geodesic
submanifolds of (real) codimension
2, but none of codimension 1 (of course, it
is far away from
being a warped product).

%\paragraph
\subsection{From the local to the global in the analytic case}
 In
presence of
(real) analyticity, a somewhere local warped product leads to
an everywhere local warped product, and then to
a global warped product in the universal cover, and finally to
  full completeness.
%the above results become global:

\begin{theorem}
\label{global}
Let $M$ be a {\bf compact} (real) {\bf analytic} Lorentz manifold.
 Suppose that
${\cal W}(M)$ is nonempty.
Then,  the universal cover
$\tilde{M}$ is isometric to a warped product of
a complete simply connected Lorentz manifold
$\tilde{N}$ of constant {\bf nonpositive} curvature and
dimension $\geq 3$, by a {\bf complete} simply connected
Riemannian
manifold $\tilde{L}$.
Furthermore,  $M$ is
(geodesically) {\bf complete}, and admits another
metric for which
$\tilde{M}$
is isometric to the direct product
$\tilde{L} \times \tilde{N}$.
\end{theorem}

The last part of this  theorem contains in particular, Carri\`ere's theorem,
and its adaptation by B. Klingler, on completeness of compact
 Lorentz manifold
%(whose proofwas adapted
of constant curvature \cite{Car} \cite{Kli}. We notice however, that we
don't reprove
Carri\`ere's theorem here, but instead use it, by  observing that its proof
may be
adapted to the general situation in the theorem above.

%\paragraph
\subsection{Warped product or local  bi-polarization, when the isometry group
is noncompact}
\label{introduction.def.2}

Let $M$ be a compact Lorentz manifold, such that
Isom$M$ is noncompact. From Part I, there exists
   at least one geodesic lightlike codimension 1
foliation of $M$.

Therefore,
for
all $x $ in $ M$, card$ C_x \geq 1$. This fact alone
may also be proved in a straightforward way, by looking  at  limits
of  graphs of
the elements of Isom$M$.

Fuschian-like behavior of Isom$M$ was described
in Part I,
having as a consequence
% leading in particular (roughly speaking) to
a dichotomy (roughly speaking):  card$C_x \leq 2$, or
$C_x$ infinite. From
the results above, the last situation
implies that ${\cal W}(M) \neq \emptyset$,
and thus $M$ has the nice structure described above,
in the analytic case.

\begin{theorem}
\label{classification}
Let $(M, g)$ be a {\bf compact} (real) analytic Lorentz manifold,
such that Isom$M$ is noncompact. Then,  exactly one of
the  two following possibilities
% (1) or (2)
 holds:

\medskip
\noindent
1)  There exists a new metric $g^\prime$
on $M$ such that:

(i)  Isom$(M, g)$ is
a normal cocompact subgroup of
Isom$(M, g^\prime)$.

%(i)  Isom$(M, g) \subset$ Isom$(M, g^\prime)$,
%and Isom$(M, g^\prime)/$ Isom$(M, g)$ is compact.

(ii)
The universal cover
%$\tilde{M}$
of $(M, g^\prime)$
is isometric to
a direct product
 $\tilde{L} \times \tilde{N}$,
 where $\tilde{L}$ is a complete simply connected Riemannian
manifold, and $\tilde{N}$ is
a complete simply connected Lorentz manifold
 of constant {\bf nonpositive} curvature and
dimension $\geq 3$.

\medskip
\noindent
2) There are two
%(may be non distinguish)
(not necessarily distinct) codimension 1
lightlike geodesic foliations ${\cal F}_1$
 and ${\cal F}_2$, such that:

(i) Any
lightlike geodesic hypersurface in $M$, is
contained in a leaf of ${\cal F}_1$
or ${\cal F}_2$. In particular, any {\bf local}
isometry of $M$ preserves each of these
 foliations, or exchanges  them.

(ii) For each $i$, there
is an analytic structure on
(the topological manifold) $M$
in respect to
which
%such that
 ${\cal F}_i$
%becomes an
 is an  analytic foliation.

\end{theorem}

Some comments are in order:\\

\paragraph{Local bi-polarization.} Recall from Part I, that $M$
is called {\bf bi-polarized} if its isometry
group
is noncompact and preserves a pair of
lightlike geodesic foliations. The situation
(2, (ii)) in the theorem above, suggests the
definition of  a local  version of this notion
(i.e.
by means of the  pseudo-group
of local isometries).
We will say that
$M$ is {\bf locally bi-polarized} if
there are two lightlike geodesic foliations ${\cal F}_1$
and ${\cal F}_2$ such that
for all $x \in M$, $C_x = \{ (T_x{\cal F}_1)^\perp,
(T_x{\cal F}_2)^\perp) \}$.

In order
%to gofurther towards
to get closer to a classification of
compact Lorentz manifolds with noncompact isometry group,
%it is surely interesting to investigate
the investigation of the geometric and
dynamical structure of
locally bi-polarized manifolds is clearly
of interest.

Let's give some  examples of locally
bi-polarized manifolds.
 Consider $SL(2, {\bf R})$
endowed with its Killing form. It has constant negative curvature,
and thus, it is by no means bi-polarized.  Now, endow $SL(2, {\bf R})$
with a left invariant Lorentz metric derived from a given
Lorentz scalar product $<,>$ on the Lie algebra
$sl(2, {\bf R})$. Suppose that $<, >$
is a kind of ``Berger's metric'', that is, it is
given by scaling
 a hyperbolic element $u \in sl(2, {\bf R})$
(i.e. $\exp tu$
is a hyperbolic one
parameter group) by a nontrivial factor, and keeping
the Killing form
on
$u^\perp$. The  metric so obtained
admits as a  local bi-polarization
the stable and unstable foliations
determined by
$\exp tu$ (see \cite{G-L} about left invariant
metrics on $SL(2, {\bf R})$).

\paragraph{The isometry group in the warped product case.} Of course,
 being locally bi-polarized
is stronger than being  bi-polarized. For
instance, ``purely'' irrational
% generic (i.e. purely irrational)
 flat structures (which
are of course not locally  bi-polarized) on
the torus, with noncompact
isometry group,  are bi-polarized (see Part I, \S 15). Also, nonhomogeneous
3-anti de Sitter manifolds are bi-polarized, whenever they
have a noncompact isometry group (this isometry group
is in fact, up to a finite index, isomorphic to
${\bf R}$, and hence it is amenable, which
implies from Part I that the underlying manifold is
 bi-polarized).

Furthermore, we
 observed in Part I, that if $M$ is not bi-polarized, and
  the factor $\tilde{N}$ in Theorem \ref{classification} has
constant negative curvature (that is,
$\tilde{N}$
is an anti de Sitter space), then,  dim$\tilde{N} = 3$
(i.e. $\tilde{N} = \widetilde{SL(2, {\bf R}) }$). Equivalently,
if dim $\tilde{N} >3$, then,  $M$
is bi-polarized. However, it seems that
this situation never happens. That is, if
$\tilde{M} = \tilde{L} \times \tilde{N}$, where
$\tilde{N}$ is an anti de Sitter space
of dimension $>3$, then Isom$M$
is compact (from our definition, we don't call
$M$
bi-polarized in this case). This was proved in the case
 that the factor $\tilde{L}$
is trivial, so  $M$ is an anti de Sitter manifold of
dimension $>3$ \cite{Zeg.anti.de.Sitter}.

%LOL REDUCTIVE HOLOLOMY

\paragraph{Regularity.} The lightlike geodesic foliations found in
Part I are a priori only Lipschitz. However, in each
of the cases, warped product or locally bi-polarized,
there are extra reasons leading to higher regularity. Indeed, in
the warped product case, we  essentially deal with {\bf global}
lightlike  geodesic foliations of the anti de Sitter or
the Minkowski spaces. They
 are easy  to
handle, and can be shown to
be analytic.

Here
%follows
is the  idea of the proof of  regularity in the
locally  bi-polarized case (which is behind the
property (2, (ii)) of Theorem \ref{classification}.

Observe that the graph (as a section) of  a codimension
1 lightlike
 geodesic foliation ${\cal F}$ on $M$,
is
%gives rise to
a Lipschitz
submanifold $P({\cal F})$,  homeomorphic to $M$, contained in
$Gr^0(M)$,
%the integrability set ${\cal I}(M)$, of
%the tautological geodesic plane field in
the Grassman bundle of lightlike hyperplanes tangent to $M$.
In fact $P({\cal F})$
is contained in the subset
${\cal D}$, the integrability domain of
the tautological plane field (see \S \ref{tautologic}),
defined by
${\cal D} = \cup_{x \in M} C_x^*$, where
$C_x^* \subset Gr^0_x(M)$
is the dual of
$C_x$ (see above). But, $\cal D$
is an {\bf analytic set}.
So, amusingly, when $\cal D$
is poor, say it equals
$P({\cal F}_1) \cup P({\cal F}_2)$, then we
win regularity for
${\cal F}_1$ and ${\cal F}_2$, because their graphs are
open in an analytic set. This implies that
${\cal F}_1$
and ${\cal F}_2$ are ``essentially''
analytic. But,   a priori,
${\cal D} $
may have
%in the general case, one can't avoid
an intrinsic singularity locus,  or a vertical locus
(where it is regular but tangent to the vertical).
However, we guess,
none of these singularities may occur in our situation, and
the foliations are actually analytic.
Anyway, we have the following
corollary of Theorem \ref{classification} (2, (ii)).

\begin{theorem} Let $M$
be a compact
%analytic
topological manifold
% admitting
which has no codimension 1 analytic foliation,
%in the sense of an
for any analytic structure on $M$.
%(for example $M$ simply connected).
Then,
any analytic Lorentz metric on
$M$ has a compact isometry group.
\end{theorem}

A classical result of A. Haefliger (see for
example \cite{Hea})
states that  compact
simply connected
  manifolds
satisfy the condition of the theorem.
Therefore, they have
compact isometry
groups.
This gives, another proof of
G. D'Ambra's theorem  \cite{D'A}, without
using Gromov's theory of rigid transformation groups.

Recently, T. Barbot \cite{Bar} has
found another examples of manifolds
satisfying
the condition of the theorem above.
For instance, a compact manifold
with a fundamental
group isomorphic to
a finite index subgroup of $SL(d, {\bf Z})$, with
$d \geq 3$, has no analytic foliation. Therefore, such a manifold
has a compact isometry group, when endowed with an analytic
Lorentz metric.

\paragraph{A properness theorem.} In fact, as  was said in
Part I, one
 may ask for a more stable compactness
of isometry groups.
% of compact simply connected
%Lorentz manifolds.
Our method
allows us to prove
 the following  result,
%generalization of D'Ambra's result,
which we state here without further details.

\begin{theorem} Let $M$ be a compact manifold, which is simply connected,
or has a fundamental
group isomorphic to
a finite index subgroup of $SL(d, {\bf Z})$, with
$d \geq 3$.
 Denote by Lor$^{\omega, 2}(M)$
(resp.  Diff$^{\omega, 2}(M)$) the space of
analytic Lorentz metrics on $M$
(resp. analytic diffeomorphisms of $M$)
endowed with the $C^2$ topology. Then
 Diff$^ {\omega, 2}(M)$ acts properly
on  Lor$ ^{\omega, 2}(M)$.

\end{theorem}

\paragraph{An application: action of Lie groups}
Due to the  works
\cite{Zim}, \cite{Gro}, \cite{D-G}, \cite{Kow},
\cite{A-S.1}, \cite{A-S.2}, \cite{Zeg.espacetemps.homogenes}
and \cite{Zeg.identity.component}, many things are now known about
 isometric actions
of connected  Lie groups on compact Lorentz manifolds. For
example, we know that, if
the affine group (of the line) $AG$
acts isometrically , then, essentially, this
action may
be extended to
$SL(2, {\bf R})$ (see \cite{A-S.2} or \cite{Zeg.identity.component}
for the correct  statement of
this fact). Let's see how to deduce  this fact,
in
the analytic case,
from
Theorem \ref{classification}.
 Of course  this  would follow by standard algebraic
manipulations, if we already knew
that the manifold was  a warped product, as described
in the point (1) of Theorem \ref{classification}. So,
%in the analytic case
%an alternative
%proof of this result
 it suffices to
show that a manifold endowed with
an $AG$-action is not bi-polarized. For
this, let $h^t$ and $T^t$ be
two one-parameter groups generating
$AG$ such that
$h^sT^th^{-s} = T^{te^s}$. Then,
two isometric flows determined by  two different hyperbolic
one-parameter groups  $T^{s_1} h^t T^{-s_1}$ and
$T^{s_2} h^t T^{-s_2}$,
 determine two different  approximately stable
foliations. Indeed the tangent bundles of these
foliations, are the orthogonal
of the negative Lyapunov spaces (which are
isotropic of dimension 1) associated to
the given flows.

%\paragraph
%\subsection{The smooth case}

\section{Geometry of warped products}
\label{warped.def}
Here, we will present  standard geometric notions related
to
warped products, and state  geometric
%criterion leading to
 criteria for the existence
of such   structures.
We will try to avoid the  use  of local calculus, and
 instead, we will use  synthetic arguments.

\paragraph{Umbilical and geodesic submanifolds.}
Let $M$ be a Lorentz manifold.
Let $S $ be a nondegenerate
submanifold of $M$, that is the metric restricted to
$T_xS$ is nondegenerate for any $x \in S$.
Recall that
$S$ is umbilical, if and only if for any
$x\in S$, the second fundamental
form $II_x$ (which is
well defined because of the non-degeneracy  hypothesis)
 has the form
$II_x = <,>n_x$, where $n_x$ is some   normal vector to
$T_xS$. The geodesic case corresponds to
$n_x= 0$,  for all $x \in S$.

Let $x \in S$,
$u \in T_xS$, and let $\gamma : ]-\epsilon , +\epsilon[
\to M$ be the geodesic in
$M$
determined by $u$. For $S$ geodesic, the image
of $\gamma$
is contained
in $S$, for $\epsilon$
sufficiently small. This fact is  true also when
$S$ is umbilical, if in addition
$u$ is {\bf isotropic} (this is a remarkable
rigidity fact in Lorentz geometry, which
has no equivalent
statement  in
Riemannian geometry) .

For  example, take $M$ to be
Minkowski space, i.e. ${\bf R}^n$
endowed with a Lorentz form $q$. The geodesic hypersurfaces
are contained in affine hyperplanes. The umbilical
hypersurfaces are contained in quadrics
$q(x-O) = c$, where $O \in {\bf R}^n$ and $c$ is
a constant (the proof is
formally the same as in the
Euclidean case).
 One can verify  that such a
quadric is  ruled, that is, it  contains the isotropic
lines which are somewhere
tangent to it.

\paragraph{Umbilical and geodesic foliations.}
A foliation is called geodesic or umbilical, if and only if its leaves
are geodesic, or umbilical. The following
is a standard fact \cite{Mol}.

%Let $M$ be a Lorentz manifold

%The following geometric
% discussion is local.
%, and so, for the sake of simplicity, we won't mind on
%topological complications

\begin{fact}
\label{geodesic.critere}
Let $\cal F$ be a foliation of $M$
such that the orthogonal $T{\cal F}^\perp$
is
integrable, that is, it  determines a foliation ${\cal F}^\perp$, say.
Then ${\cal F}$ is geodesic (resp. umbilical)
if and only if the holonomy maps of the foliation ${\cal F}^\perp$,
seen as local diffeomorphisms between leaves of ${\cal F}$,
preserve the metric (resp. the conformal structure)
induced on these leaves (of ${\cal F}$).

Let $F$ and $F^\perp$ be the local leaves
through  some point
for $\cal F$
and
${\cal F}^\perp$,
respectively, leading
to a local diffeomorphism of $M$
with $F \times F^\perp$. If
$\cal F$
is
geodesic, then the metric has the form $m_{x, y}=
h_x \bigoplus g_{(x, y)}$, where $h$ is a metric
on $F$ and for any $x \in F$,
$g_{(x, .)}$ is a metric on $F^\perp$.

If ${\cal F}$ is umbilical,  then the
metric has the form $m_{x, y} =
w(x, y)h_x \bigoplus g_{(x, y)}$,
where $w $ is a
function on $F \times F^\perp$.
\end{fact}

 \paragraph{Warped products.} As  was said in the introduction,
 a Lorentz manifold $M$
is  a warped product
of a Lorentz manifold
$(N, g)$ by a Riemannian manifold
$(L, h)$, if $M$ is isometric to
the product $L \times N$, endowed with
a metric of the form
$h \bigoplus w g$, where $w$ is
some positive function defined
%(only)
 on $L$.
 We call $N$ the {\bf normal} factor of the warped
product.

The factors $L$ and $N$ define two foliations
denoted
by
${\cal L}$
and ${\cal N}$
respectively;  ${\cal N}$ is called
the {\bf normal foliation}
of the warped product.
From the form
of the metric and the above discussion, we
infer
%that
%have:
%These two foliations satisfy
the following geometric properties:
${\cal
L}$ is geodesic, and ${\cal N}$
is umbilical
%(by this we mean that their leaves
%are geodesic and umbilical respectively).

In terms of holonomy pseudogroups,  this means that
 the holonomy of ${\cal N}$ (resp. ${\cal L}$)
 preserves the transverse metric (resp. the transverse
conformal structure).
 In fact, to
characterize warped products, we just need that the holonomy
maps of $\cal L$ have  constant conformal
distortion, that is,  they
are {\bf homothetic}. In particular the holonomy of ${\cal L}$
is projective, i.e. maps geodesics to geodesics.
Here is  a related stronger
property:
%other words, let

\begin{fact}
\label{geodesic.submanifold}
A submanifold $S $ of
$N$ is  geodesic
in  $N$, if and only if  $L \times S$
is geodesic in $M = L \times N$.

\end{fact}

\begin{proof} It suffices  to consider
the case dim$S = 1$.

Suppose that
$S$
 is geodesic in $N$. By considering the family of normal
geodesic of a hypersurface orthogonal to $S$, we
can locally extend $S$ to a geodesic foliation
${\cal S}$ of $N$, admitting an orthogonal foliation
${\cal S}^\perp$.
 Consider the foliation ${\cal F}$
of $M$ with leaves of the form $L \times {\cal S}_y$
where ${\cal S}_y$ is a leaf of ${\cal S}$. It
has a normal foliation ${\cal F}^\perp$
with leaves $ \{x\} \times {\cal S}^\perp_y$.
It then follows that a holonomy map
of ${\cal F}^\perp$
has the form $\psi: (x, y) \in L \times {\cal S}_{y_1} \to
(x,  \phi(y)) \in L \times {\cal S}_{y_2}$,
where $\phi$ is a holonomy map
of ${\cal S}^\perp$, which is isometric by
hypothesis. Because the metric
on $\{x\} \times N$
is a constant times the metric
of $N$, $\phi$ sends
$\{x\} \times {\cal S}_{y_1}$ isometrically
onto $\{x\} \times
{\cal S}_{y_2}$. Therefore $\phi $
is isometric. Hence from Fact \ref{geodesic.critere}, ${\cal F}$
is geodesic, and in particular
$L \times S$ is geodesic in $M$.

For the converse, that is,  $L \times S$ being geodesic
implies that $S$ is geodesic, we use the following
 general fact. Its proof follows from a
standard calculation.

\end{proof}

\begin{fact} Let $A$ and $B$
be  submanifolds of $M$, and suppose that
$B$
is geodesic.
 Suppose that $B$ is transverse and orthogonal to
$A$, that is, for all
$x \in A \cap B$, $T_xB$ contains $(T_xA)^\perp$.
 Then $B\cap A$ is geodesic in $A$.
\end{fact}

%Here we are interseted in
%the case $N$ is a Lorentz manifold
%with constant curvature.

\paragraph{Criterion for warped products.} Here is the proposition which
we will apply to prove the existence of warped products
in \S \ref{tautologic}.

\begin{proposition}
\label{warped.critere}
Let $M = L \times N$ be endowed with a metric
such that the foliation $\cal L$
(resp. $\cal N$) is geodesic (resp. umbilical).

Suppose that for
all $(x, y) \in L \times N$, there are  geodesic hypersurfaces in $M$,
$H^1,
\ldots, H^d$,
containing $(x, y)$ and such that:

i) $H^i$ is invariant by the foliation $\cal L$
(i.e. it is a union of leaves of $\cal L$), and

ii) The directions  $(T_{(x,y)}H^1)^\perp \cap T_yN , \ldots,
(T_{(x,y)}H^d)^\perp \cap T_yN$
generate $T_yN$.

Then the leaves  $\{x\} \times
 N$ have constant curvature, and
$M$ is a warped product.

\end{proposition}

\begin{proof}
 One can write $H^i= L \times S^i$, where
$S^i$ is a hypersurface of $N$.  From Fact \ref{geodesic.submanifold},
$\{z\} \times S^i$  is a geodesic hypersurface
in  $\{z\} \times N$, for all $z \in L$.

Thus, $\{ z\} \times N$, admits, at each
point, geodesic  hypersurfaces, whose orthogonal
directions generate the tangent space at each point.
It will be shown at Proposition \ref{constant.critere},
that this implies that
 $\{ z\} \times N$ has constant curvature.
  We may assume that
the sign
of
 curvature
 is independent of $z \in L$.

A holonomy map (of $\cal L$)
taking  $\{x\} \times N $ to $ \{x^\prime\} \times N$,
maps $\{x\} \times S^i  $ to $ \{x^\prime\} \times S^i$.
One may then call it ``partially projective''.

The warped product
property means that any  holonomy map (of
${\cal L}$)
%, defined between
%leaves of the form
$\{x_1\} \times N \to \{x_2\} \times N$,  is  homothetic.
By hypothesis, $\cal N$ is
umbilical, and hence
these holonomy maps are conformal.
The question then becomes, is a  conformal and
``partially projective'' map between two Lorentz
manifolds,  homothetic?
  In our case, the leaves
$\{z\} \times N$
 have constant curvature of the
same sign, and therefore are  homothetic.
 The proof of
 the warped product
property  can  then  be achieved with help of the following fact.
\end{proof}

\begin{fact}
\label{homothetic.critere}
 Let $N$ be a Lorentz manifold of constant curvature,
with dimension $\geq 3$. Let
 $\phi: N \to N$ be a conformal local diffeomorphism, which satisfies
the following condition.
For
any  $y \in N$, there is $S^1 \ldots,
S^d$,  geodesic hypersurfaces
 containing $y$,
such that $(T_xS^1)^\perp, \ldots, (T_xS^d)^\perp$
generate $T_xN$, and such that,
 their images $\phi(S^i)$ are geodesic.
 Then $\phi $ is a homothety. (In fact
$\phi$ is an isometry unless $M$ is flat). (Interpretation:
a conformal and ``partially projective'' transformation
is  homothetic).

\end{fact}

\begin{proof}
We think that the interpretation
of the fact is sufficiently
convincing,  and so,
%for the sake of simplicit
to avoid complicated notation, we restrict
ourselves  to the flat case. Also, we
will assume that the involved geodesic
 hypersurfaces $S^1, \ldots, S^d$
are lightlike, because,  that is what we need
for application, in the present paper.
  Thus $M$
is
the Minkowski space
${\bf R}^{1, n-1}$. By composing
with an isometry,
we may
suppose that $\phi $ is a local
conformal
diffeomorphism fixing $0$, and
that $D_0\phi$
is a homothety. In particular
$D_0\phi$ keeps invariant
any tangent line at $0$.
It then follows that $\phi$
keeps
invariant each isotropic line through $0$, since
conformal diffeomorphisms preserve isotopic geodesics.
Furthermore, $\phi$ keeps invariant the
geodesic hypersurfaces $S^1, \ldots, S^d$.
(because they are sent by $\phi$
to geodesic hypersurfaces with the same tangent space).
Also, by considering intersection of sub-families
 of these hypersurfaces, we
infer the existence of a basis
$\{v_1, \ldots, v_n \}$ of $T_0{\bf R}^n$,
determining lines kept invariant by $\phi$.
These vectors are spacelike, that
is $<v_i, v_i> >0$,  since
the hypersurfaces $S^i$ are lightlike.

Let $P$ be
an affine 2-plane, so  $\phi(P)$
is 2-dimensional generalized sphere, i.e. an
affine plane, or
 a quadric defined by means
of the Lorentz form (this
% because
 follows from the fact  that $\phi$ is conformal,
as in the Riemannian case).

There are two possibilities
for a generalized 2-sphere which is not an affine plane. If
it is (somewhere and hence everywhere) timelike (i.e.
the induced metric on it, is of Lorentz
type,  then it is ruled, by means of
 a pair of foliations by   isotropic lines.
In contrast, if it is spacelike (i.e.
the induced metric  on it is Riemannian), then it
 contains no line.

Let $P$ be a spacelike affine plane which
contains a line defined by some $v_i$. Then
$\phi(P)$
is a  spacelike generalized 2-sphere
which contains a line. From the above
discussion, $\phi(P)$ must be an affine plane.
But, since $D_0\phi$ is a homothety, we have
$\phi(P) = P$.
Thus all the affine spacelike 2 planes containing
a line determined by some $v_i$ are invariant by $\phi$.
A standard analyticity argument
show that this extends to all the affine 2-planes
 without the spacelike condition.
Taking the  intersection of these planes, for
various $v_i$, and again by an analyticity argument,
we  conclude that every line through
$0$ is invariant by $\phi$.

But the same argument works for any point of $M$ (by
composing with an appropriate
isometry). This means that $\phi$
is projective. In particular the restriction of $\phi$
to an affine plane $P$ through $0$ is conformal and
projective. This is equivalent to
a conformal and projective local transformation of the
Euclidean plane ${\bf R}^2$, if $P$ is
spacelike. This is easily seen to be a homothety.
Because, there are many
such planes, one  concludes that $\phi $ itself
is a homothety.

\end{proof}

\section{The tautological geodesic plane field}
\label{tautologic}

%\subsection{General setting}
%\paragraph{General setting.}

An affine connection (e.g.
a pseudo-Riemannian structure) on a manifold $M$,
permit to define a
tautological
geodesic plane fields
on the Grassmann bundles
of tangent $k$-planes
$Gr_k(M) \to
M$.
This generalizes the classical
construction
of the geodesic flow   for $k=1$.
The connection yields a horizontal bundle $H$, supplementary to
the vertical.  For $p \in (Gr_k)_x$, we identify
$H_p$ with $T_xM$. Then ${\cal P}_p \subset H_p$
 is identified with $p \subset T_xM$.
 In general, ${\cal P}$ is not integrable. Indeed
one may prove (see below), for $M$ pseudo-Riemannian,  that if
${\cal P}$ is
integrable for $k \neq1$ and $k\neq $ dim $M$, then
$M$ has constant curvature.  (It seems that when $M$
is merely affine, then the conclusion
is that  $M$  is projectively flat).

Observe that (maximal)
integral submanifolds
of ${\cal P}$ project on geodesic
submanifolds of dimension $k$ in $M$.
 Conversely, if $L$ is a geodesic $k$-submanifold
of $M$, then  the image of
the Gauss map $x \to T_xL$ is
 an integral submanifold of ${\cal P}$.

We denote by $\exp: TM \to M$,
the exponential map defined on its domain of
definition,  which is a neighborhood of the
 zero section. For $p \in Gr_k (M)$, let $\exp p$ denote
 the
image by $\exp$ of the intersection
 of $p$
with the domain of definition of $\exp$.

\begin{definition} The domain of integrability
of ${\cal P}$ is the set of $p \in Gr_k(M)$
such that, if $p \subset T_xM$, then
a neighborhood of $x$ in  $\exp_x p$ is geodesic.
%It is a closed subset of
%$Gr_k(M)$, on which ${\cal P}$
%induces a lamination.
\end{definition}

\subsection{The Grassmannian of lightlike hyperplanes
of a Lorentz manifold}
Let $M$ be a Lorentz manifold
and denote by $Gr^0(M) \subset Gr_{n-1}(M)$ the Grassmannian
of {\bf lightlike hyperplanes}  tangent to
$M$ ($n $ is the dimension of $M$).

The tautological geodesic plane field on
$Gr_{n-1}(M)$
 is tangent to $Gr^0(M)$
We will denote the
tautological plane field restricted to
$Gr^0(M)$ by $\tau$, and
by ${\cal D}$ its
 {\bf integrability domain}.

Therefore, as in \S \ref{introduction.def.2},
  the fiber  ${\cal D}_x $, is the dual of $C_x$.
 Recall that a direction $u $
 belongs to
$C_x$, if and only if $u^\perp$ is tangent to
a lightlike geodesic hypersurface.

 To start with, notice the following
rigidity.

\begin{proposition}
\label{constant.critere}
 Let $M$ be a Lorentz manifold
of dimension $\geq 3$,
 such that $C_x$ generates $T_xM$,  for any $x \in M$.
 Then $M$ has constant curvature.

\end{proposition}

\begin{proof}
The hypothesis means that
 for any  $x \in M$, there
are  $H^1,  \ldots,
H^d$,  {\bf geodesic lightlike hypersurfaces}
 containing $x$,
such that $(T_xH^1)^\perp, \ldots, (T_xH^d)^\perp$
generate $T_xM$.

Fix $x \in M$, denote
$T_xH^i$ by $B^i$, and choose $b^i$ an isotropic vector
such that $B^i = (b^i)^\perp$.

 For  $u \in T_xM$,
denote by $A_u$ the curvature operator $A_u: v \in
T_xM \to R(u, v)u \in T_xM$. Then,
for $u  \in B^i$, $A_u$
preserves $B^i$ (since  geodesic
submanifolds are ``invariant'' by the curvature operator).
Moreover, $A_u (b^i)$ is collinear
to $b^i$ (for $u \in B^i$).
Indeed $<A_u(b^i), v>= <A_u(v), b^i>$. The last
quantity equals $0$ if $v \in B^i$ since
$A_u(v) \in B^i$, and hence
$A_u(b^i) \in (B^i)^\perp
 = {\bf R} b^i$.

Choose $e_i$ a unit director vector of $\cap_{j\neq i}B^j$, and
consider $A_{e_i}$.  Since $e_i \in B^j$, for
$j \neq i$, there is $\lambda_{i, j}$,
such that
 $A_{e_i} b^j = \lambda_{i, j} b^j$
(for $i \neq j$).

Since $A_{e_i}$
is symmetric, $ \lambda_{i,j}<b^j,b^k> =
<A_{e_i}b^j, b^k> =
 <A_{e_i}b^k, b^j>
= \lambda_{i,k}<b^j,b^k>$, and we
have $
\lambda_{i,j} =
\lambda_{i,k}$
 ( for
$j \neq k$,
$<b^j,b^k>  \neq 0$, because
both $b^j$
and $b^k$
are isotropic). Write $\lambda_i = \lambda_{i, j}$.
Thus, the sectional curvature of any nondegenerate plane
which contains $e_i$ equals $\lambda_i$. From this, we infer
that $\lambda_1= \ldots= \lambda_n $
% $ = $, say,  $\lambda$
(to see this,  consider
2-planes generated by two vectors $e_i$ and $e_j$).
 One may  use standard algebraic manipulations
to show that all the 2-planes in $T_xM$ have the same
sectional curvature, and then deduce from
Schur's lemma that $M$ has constant curvature.

\end{proof}

\subsection{Main result}

\begin{theorem}
\label{hyper.warped}
 For $x \in M$, denote by
$E_x$ the linear space generated by
$C_x$.
Suppose that
 for  an open subset $U \subset M$, we have:

i) $x \in U \to  E_x$
%(the vector space generated by $C_x)$)
 determines  a smooth plane field of dimension $\geq 3$, and

ii) card$C_x > $ dim$E_x$, for $x$ in a dense subset of $U$.

 Then $E$
determines a local  warped
product structure on
$U$ (i.e. $E$ is tangent to
the normal foliation
of a local warped structure on
$U$), with the leaves of $E$
having constant curvature.

% Then $E$ and
%$E^\perp$
%determine a local  warped
%product structure on
%$U$, with the leaves of $E$,
%having constant curvature.

\end{theorem}

\begin{proof}
There are several  steps.

\begin{step} $E^\perp$
is integrable and has geodesic leaves. We denote its
tangent foliation
by
$\cal L$.

\end{step}

\begin{proof}
Let $H^1, \ldots H^d$, be
geodesic lightlike hypersurfaces containing a point
$x \in U$, such that
$(T_xH^1)^\perp, \ldots, (T_xH^d)^\perp$ generate
$E_x$. Denote  $L = \cap_{i}H^i$, and let's show that it
is a leaf of $E^\perp$.

 Denote by $X^i$ a nonsingular isotropic vector
field tangent to $H^i$ ($X^i$ is defined along
$H^i$). Then for $y \in H^i$,
$X^i(y) \in C_y$.  Hence, if
$y \in L = H^1 \cap \ldots H^d$, then
, $X^1(y), \ldots, X^d(y) \in C_y$. Thus, by continuity of
$E$, $E_y$ is generated by $X^1(y), \ldots, X^d(y)$.
Observe now that $T_y L =
\cap_i T_yH^i = \cap_i (X^i(y)^\perp =
(\Sigma_i {\bf R}X^i(y))^\perp = E_y$. Therefore $L$
is a leaf of $E^\perp$ containing $x$, which   (being
 an intersection of geodesic hypersurfaces) is a geodesic
submanifold.
\end{proof}

%\paragraph{Weingarten's endormophisms for plane fields.}
\begin{step}
Weingarten's endomorphism for plane fields.
\end{step}
For a plane field  $x \to G_x \subset T_xM$, such
that the metric restricted to
$G_x$ is {\bf nondegenerate},  one defines a second fundamental
form
and Weingarten's endomorphism as follows.
For $X$ and $Y$ vector fields
tangent to $G$, and $Z$ a vector field orthogonal to
$G$,
$II: G \times G \to G^\perp$, and
$A_Z: G \to G$,
are defined by
the equalities:  $<II(X, Y), Z> = <\nabla_XY, Z> = <-A_ZX, Y>$.
We have: $0= X<Y,Z>= <\nabla_X Y, Z> + <Y, \nabla_XZ>$, and
hence $A_Z(X) $ is just the projection of
$\nabla_X Z$ on $G$.
It turns out that $II$ and $A_Z$ are tensorial, that is,
they depend only on the pointwise values of $X$, $Y$ and
$Z$.  Notice  the following property:

\begin{fact} A plane field $G$ is integrable, if and only if its second
fundamental form, or equivalently its Weingarten's endomorphisms,
are symmetric.

In particular, if any  Weingarten's endomorphism
of $G$ is a {\bf homothety} (that is,  it induces
a scalar multiplication on $G$), then,
$G$ is integrable, and has umbilical leaves.

\end{fact}

Our plane field $E$ is nondegenerate, since it contains at
least two isotropic directions.

%\paragraph{End of proof of the theorem, assuming that
%all Weingarten's endomorphism are homotheties.}

\begin{step}
End of proof of the theorem, assuming that
all Weingarten's endomorphism are homotheties.
\end{step}

\begin{proof}
In this case,
from the fact above, $E$ is integrable and has umbilical leaves.
Then we use Proposition \ref{warped.critere}, to
deduce that $E$ and $E^\perp$ give rise to a warped product.
Indeed, $E$ is umbilical,  $E^\perp$ is geodesic,
and furthermore, the condition on the existence of
geodesic hypersurfaces saturated by $\cal L$, is
well satisfied. Indeed,  the hypersurfaces,
$H^1, \ldots, H^d$, introduced in the beginning of the proof
of the integrability of $E^\perp$, are saturated by $E^\perp$.
\end{proof}

%the partial projectivity'' condition on the
%holonomy of
%$E^\perp$ is well satisfayed.

\medskip

 Now follow the steps
of the proof that
the weingarten's endomorphisms are actually homothetic.

We will consider eigenspace
splittings   for $E$, and then splittings of the factors of the initial
splitting, and so...
All these splitting are smooth, if
 we  restrict ourselves to  an open dense subset.

 Notice that
this doesn't   loss of  generality.
 Indeed,
if we are able, at the
final stage,  to prove that the weingarten's endomorphisms
 are homothetic, in a dense set,
%at generic point,
then they will
be homothetic everywhere. So, in the  sequel, we will
always suppose that we
are near a generic point.

%\begin{fact}

\begin{step} \label{step.1} Let  $Z$
be a smooth vector field tangent
to
 $E^\perp$($=T{\cal L}$).  Then, for any
$u \in C_x$, $u^\perp \cap E_x$ is invariant by
the weingarten's endomorphism $A_Z(x)$ (or
equivalently, $u$
is an eigenvector of
the dual weingarten's endomorphism $A_Z^*(x)$).

\end{step}

\begin{proof}
Let $u \in C_x$, and consider $H$ a
lightlike geodesic hypersurface
such that $T_xH = u^\perp$.
Observe that $Z$ is tangent to
$H$
(over points of $H$).  This is because
$E ^\perp$ itself is tangent to $H$ (or
equivalently $H$ is saturated by the foliation
${\cal L}$). Since $H$ is geodesic,
%we have:
the covariant derivative $\nabla_X Z(x)$
belongs to  $ T_xH$, for $X \in T_xH = u^\perp \cap
E_x$, and hence, its
 projection $A_Z X(x)$  belongs to
$u^\perp \cap E_x$. That is,
$A_Z(u^\perp \cap E_x) \subset u^\perp \cap E_x$.

\end{proof}

\begin{step}
\label{step5} The 3-dimensional case.
\end{step}
To start with,  let's give the proof when
 dim$E= 3$. The cardinality condition
in the theorem  means that card$C_x \geq 4$, for
 any $x \in U$. Hence $A_Z(x)$ has at least
4 {\bf isotropic} eigenvectors. Thus, since
dim $E =3$, the eigenspace decomposition
of $E_x$ has at most two factors $E_x = A \bigoplus B$.
If this is nontrivial, then up to a switch of factors, we have
dim$A= 2$ and dim$B=1$. Of course our isotropic eigenvectors
belong to
$A \cup B$. But, for
dimensional reasons,  $A$ contains at most
 2 isotropic directions, and
$A $ contains  at most 1. This contradiction
implies that the decomposition
is trivial, that is
$A_Z(x)^*$
is a homothety, and thus, also is
$A_Z(x)$.

\begin{step} Getting a
partial
% part of $C_x$ which determines a
 warped product structure.

\end{step}

From Step
\ref{step.1}, we infer  that all the dual endomorphisms
$A_Z(x)^*$, for
$Z \in E^\perp_x$ are simultaneously diagonalizable. This
determines a splitting $E_x = E_x^1 \bigoplus \ldots
\bigoplus E_x^k$ of
common eigenspaces of all the $A_Z(x)^*$.
Denote
$C_x^i = C_x \cap E_x^i$, then,
we have:
$C_x = C_x^1\cup \ldots \cup C_x^k$,
%$C_x  = (C_x \cap E_x^1)\cup \ldots \cup (C_x\cap E_x^k)$,
since the elements
of $C_x$ are eingenvectors of $A_Z(x)^*$. Therefore,
for some factor, say, $E_x^1$, we have
card$C_x^i>$ dim$E_x^1$.  In particular, as observed in Step
\ref{step5},
dim$E_x^1 \geq 3$.

Observe now that if $H$ is a lightlike geodesic
hypersurface containing $x$,  such that $u = T_xH^\perp$ belongs
to $C_x^i$, then for
all $y$ near $x$, $T_yH^\perp$ belongs
to$C_y^i$. Indeed, $T_yH^\perp$ belongs to
$C_y =  C_x^1\cup \ldots \cup C_x^k$, and by
continuity,
%it remains in
$T_yH^\perp \in C_x^i$.

This observation allows us to prove the same properties
 for $E^1$, as this was
already proved for
$E$ itself (that is ${E^1}^\perp$ is integrable and geodesic, and
the dual weingarten's endomorphisms of $E^1$ admit
the elements of $C^1$ as eigenvectors). Also, by the same argument,
the pair $(E^1, (E^1)^\perp)$ would
determine
a warped product structure, if the weingarten's endomorphisms of
$E^1$ are homothetic. If not we get in a similar way, a
splitting of $E^1$. By induction, we arrive to a
sub-bundle $G$ of $E$, with dim$G\geq 3$, which
gives rise to
a warped product structure.

We may sum all the intermediate decompositions, and write
$E = G \bigoplus R$, where $R$ is a sub-bundle
of $E$, such that $C_x = (C_x \cap G_x) \cup (C_x \cap R_x)$.

\begin{step} Contradiction
\end{step}

As we  said in the beginning of this paper, a  (local) isometry,
%for a warped product,
of the
of the normal factor, extends to a (local)  isometry of
the warped product
 (see \S \ref{extension}).
In our case, a leaf of $G$ has constant curvature. In
particular, for any $x \in U$,  its local isotropy
group
contains
$O(1, d-1)$, where $d= $ dim $G$.

The infinitesimal action of $O(1, d-1)$
on $T_xM$  preserves $G_x$ and $R_x$.  In
fact, from the true definition of the extension of the action of
$O(1, d-1)$, its action on $T_xM$ is conjugate to
that  on ${\bf R}^{1, d-1} \bigoplus {\bf R}^{n-d-1}$,
where it acts as usual on the first
factor, and trivially on the second one. Here
$G_x$ corresponds to ${\bf R}^{1, d-1}$ and
$G_x^\perp$
to ${\bf R}^{n-d-1}$. Observe that the
subspace of fixed vectors of this action
is exactly ${\bf R}^{n-d-1}$.

To $R_x$ corresponds a $O(1, d-1)$-invariant
subspace $A$
%supplementary to
intersecting ${\bf R}^{1, d-1}$
trivially. This space
is not spacelike, since it is generated by  isotropic
vectors.
%elements ($R_x \cap C_x$ is
%supposed to be nontrivial).

To conclude, we just note that this is impossible.
Indeed, if dim$A =1$, then,  $O(1, d-1)$
acts trivially on it, since $O(1, d-1)$
 has no nontrivial  1-dimensional
representation. If not, $A$ is Lorentzian, and
thus $A^\perp$ is spacelike and $O(1, d-1)$-invariant.
But $O(1, d-1)$, as a simple noncompact Lie group,
has no nontrivial representation preserving a positive scalar product.
Thus $O(1, d-1)$ acts trivially on $A^\perp$. This contradicts
the fact that the space of fixed elements of
the representation
is exactly ${\bf R}^{n-d-1}$.

\end{proof}

\begin{remark} {\em  Although it was crucial in our proof
(especially at Step \ref{step5}), the cardinality condition
$\hbox{card}C_x > \hbox{dim} E_x$  might
perhaps be relaxed to the more natural
condition $\hbox{dim}E_x \geq 3$ (of course, by definition
we always have $\hbox{card}C_x \geq \hbox{dim} E_x$). In
fact, this is exactly the content of Proposition
\ref{constant.critere}, in the extremal case when
$\hbox{dim}E_x= \hbox{dim}M$.}
\end{remark}

\subsection{The structure of ${\cal W}(M)$.}

 Recall from \S \ref{introduction.def.1} that
${\cal W}(M)$ denotes  the open set points of $M$,
having a neighborhood isometric to a warped product
of a Lorentz manifold of constant curvature and dimension
$\geq 3$, by some Riemannian manifold.
 Observe that, a priori, a Lorentz
manifold may be written as a warped
product in
many fashions. However, the structure given by
the theorem above is unique, because it is
 associated  to the map $x \to C_x$. It
is in a natural sense the maximal
warped product structure (
among those with
a normal factor
  of
constant curvature) on $M$.

Define

$M_{smooth}= \{ x \in M /$  there is a neighborhood  $V$ of
$ x$  such that $ y \in V \to E_y  $ is smooth $ \}$, and
$M_{\leq k} =\{ x \in M/
\hbox{card} C_{x \leq k} \}$.

Note that,  if
$x \notin M_{\leq \hbox{dim} M}$, then
card $C_x >$ dim $M \geq $ dim $E_x$, and
hence the cardinality condition
of Theorem \ref{hyper.warped} is satisfied.

In addition,
$M- \hbox{int} (M_{\leq \hbox{dim} M})$
contains a dense set  of points
where card$C_x >$ dim$M$ (here int denotes
the interior).
Therefore, we have the following corollary:

\begin{theorem}
\label{hyper.warped.2}
%${\cal W}(M)$ contains
%the interior of
 int($M_{smooth} -\hbox{int}(M_{\leq \hbox{dim}M})) \subset
{\cal W}(M)$.
%In fact,  in the interior of

More precisely,
int($M_{smooth} -\hbox{int}(M_{\leq \hbox{dim}M}))$
%int $M_{smooth} -M_{\leq \hbox{dim} M}$,
 has a canonical  pair $({\cal N}, {\cal L})$
of foliations which determines a
local warped product with the leaves
of ${\cal N}$ having constant curvature. Furthermore, this pair
of foliations is invariant under  the pseudo-group
of local isometries  of $M$.

\end{theorem}

\paragraph{The analytic case.}  In this case the integrability
domain  ${\cal D}$ is an analytic set.
 The smoothness (in fact the analyticity) condition
on $E$, is always satisfied, away from some analytic sets.
Indeed, the assignment $x \to {\cal  D}_x = C_x^*$,
is analytic, in an obvious sense, away from
some analytic set.
Therefore, we have  the following corollary.

\begin{corollary}
\label{hyper.warped.analytic}
Suppose that $M$ is analytic. Then
${\cal  W}(M) \neq \emptyset$, whenever
int$(M_{\leq \hbox{dim} M})$ is not dense.
%More
%precisely, there is an open set $W$, which is
%dense in
%$M-M_{\leq \hbox{dim} M}$, in which $M$ admits a canonical
%local warped product structure.

\end{corollary}

%The same is true for the sets $C_x$.

\section{Completeness properties. Proof of Theorem \protect\ref{global}}

Here, we prove Theorem \ref{global}, which  is essentially
that if $M$ is analytic and
${\cal W}(M) \neq \emptyset$
% \Longrightarrow
then in fact ${\cal W}(M) = M$.

\subsection{Extension  of the warped product structure}
\label{extension}
 Here are our  two fundamental extension tools.

1) The first one, mentioned
in the introduction of this article,  is the extension to
a warped product of the isometries of
its normal factor. Indeed, let
 $M = L \times_wN$,
then, any (local) isometry
$f: N \to N$ induces
a (local)  isometry
% $\bar{f}: M \to M$, defined by
 $\bar{f}: (x, y) \in L \times N \to
(x, f(y)) \in L \times N$. The fact
that $\bar{f}$
is an isometry, follows from the fact that
 the metric
of $M$ has the form
 $h \bigoplus w g$, where $w = w(x)$ is
a positive function defined (only) on $L$. Thus
$\bar{f}$
preserves and acts isometrically on the leaves
of the foliations determined by  each of the factors $N$ and
$L$.

By the same rule, Killing fields of $N$ determine
Killing fields on $M$.

\medskip
2)
The second key  extension fact
is that a Killing field defined
on an open subset of  a {\bf
simply connected
analytic} Lorentz manifold, extends (as a Killing field)
to the whole manifold  (\cite{Nom}
and \cite{D-G}).  \\

Now let  $M$ be a compact analytic  Lorentz manifold
such that ${\cal W}(M) \neq \emptyset$.
Let $x_0$ a point of ${\cal W}(\tilde{M})$,
for which the dimension of the normal factor
 (with constant curvature) is  maximal (among all points
of ${\cal W}(\tilde{M})$).
Denote this dimension  by  $d$, and
let $\tilde{N}$ be the complete simply connected constant curvature
Lorentz manifold of dimension $d$, and
%locally isometric
%(i.e.
having the same scalar curvature as  the leaf of $x_0$
(in the local warped product).

Denote by $\cal G$ the Lie algebra of Killing fields of $\tilde{N}$.
From the extension facts recalled  above, there is a faithful
action of $\cal G$
on $\tilde{M}$, that is a monomorphism $X \in {\cal G} \to
\bar{X} \in {\cal K}$ where $ {\cal K}$
is
the Lie algebra of
Killing fields on $\tilde{M}$.  We denote the
${\cal G}$-orbit of a point $x$ by ${\cal G}x$.

By definition, near $x_0$, the ${\cal G}$-orbits
determine a local warped product. The goal is to prove that
the ${\cal G}$-orbits determine
a local warped product everywhere. That is, firstly,
the ${\cal G}$-action
gives rise to a regular foliation (i.e. of constant dimension),
 with leaves locally homothetic to
$\tilde{ N}$.  Secondly, the orthogonal
is integrable, and form together with the ${\cal G}$-foliation
a local warped product.
%Notice the
The analyticity reduces the proof   to the following
non-degeneracy fact.

\begin{fact}
\label{lorentzian}
Let $ U $ be
 the set of points of $\tilde{M}$
having a Lorentzian (also called timelike)
orbit, i.e. the induced metric on these orbits is of
Lorentzian type. Then $U$ is open, and
the ${\cal G}$-action
determines  a local warped product.
In particular, in order to prove that the ${\cal G}$-action
determines everywhere a local warped product, it suffices to
prove the equality: $U = \tilde{M}$.
\end{fact}

\begin{proof}
Observe that dim${\cal G}x \leq d$, for all
$x \in \tilde{M}$, with equality in
an open dense set.  This follows by analyticity (indeed,
if $X_1, \ldots, X_{d+1} \in {\cal G}$, then
$\bar{X}(x) \wedge \ldots \wedge \bar{X}_{d+1}(x) = 0$
in an open set).

Recall that, if  a pseudo-Riemannian manifold of dimension
$\leq d$, has a Killing algebra of the same dimension
as that of a manifold of constant
curvature and dimension
$d$, then this manifold has dimension $d$ and is of
constant curvature, of the same sign.

This shows in particular
that $U$ is open, and that the ${\cal G}$-orbit of
any point of
$U$
is locally homothetic to $\tilde{N}$. The orthogonal plane field
%distribution
$x \to T{\cal G}^\perp$ is analytic on $U$.  Consider
its second fundamental form:
$II: T{\cal G}^\perp \times T{\cal G}^\perp \to
T{\cal G}$. Its vanishing means that
the orthogonal is integrable and has geodesic
leaves. To check that $II_x= 0$, for
$x \in U$, we just use its equivariance
under
the action of the isotropy
algebra
$o(1, d-1)$ (in ${\cal G}$) of $x$. Indeed,
$o(1, d-1)$ acts trivially on
$T_x{\cal G}^\perp$, since it preserves
a positive scalar product on it. Therefore,
 for all $u$, $v$ $\in T_x{\cal G}^\perp$,
$II_x(u, v)$ is invariant
under  the action of
$o(1, d-1)$ on
$T_x{\cal G}$, and hence
$II_x = 0$.

Finally, to verify the warped product condition,  that
is, that any  holonomy map of ${\cal G}^\perp$ seen as
a local diffeomorphisms between two leaves of
${\cal G}$ is  homothetic, we just observe that this
holonomy map commutes with the action
of ${\cal G}$ on these two leaves.

\end{proof}

%We  start by
% considering some properties of the Killing
%algebra
%of constant curvature manifolds.

\subsubsection{The nonpositively curved case}
We have the
following  stronger result in the nonpositively curved case:

\begin{proposition}
\label{non.positive.case}
Let the Lie algebra
${\cal G} = {\cal G}(\tilde{N})$,
where  $\tilde{N}$ is  a
complete simply connected nonpositively curved
manifold, act isometrically on an analytic Lorentz manifold $
\tilde{M}$, and determine somewhere
 a local warped product, with a normal factor
locally homothetic to $\tilde{N}$. Then ${\cal G}$
determines everywhere a local warped product
($\tilde{M}$ is assumed   to be connected).
\end{proposition}

Behind this fact is the existence
in the nonpositively curved case, of lightlike Killing fields.
They don't exist at all in the positively curved case.

The proposition
 itself
is false in this case (see below).  We will prove it
(the proposition) assuming  in addition that $\tilde{M}$
is the universal cover of a compact manifold, and that
it is the action related to ${\cal W}(M)$. But
immediately after that  proof, we will show that
this situation never occurs, because such a compact manifold
doesn't exist.

\paragraph{Lightlike Killing fields.} A vector
field $V$ on a Lorentz manifold is lightlike (or isotropic)
if for all $x $, $<V(x), V(x)> = 0$.

Killing lightlike vector fields have the following
remarkable property.

\begin{proposition}
\label{non.singular}
(\cite{BEM}, \cite{A-S.1}). Let $V$ be
a nontrivial lightlike Killing field. Then
$V$ has no singularity. Furthermore, $V$ has geodesic orbits.
(In fact, more generally, $V$
is singularity free, when
it is   nonspacelike, i.e. $<V(x), V(x)> \leq 0$).
\end{proposition}

\begin{proof}  For the Minkowski space, a Killing
vector field vanishing at 0 is a linear Killing
vector field. It is tangent to ``pseudo-sphere'' $q = $ constant, where
$q$ is the Lorentz form. But for a
negative constant, this  level is spacelike (it
is a  Riemannian hyperbolic space). Hence, the Killing field
is somewhere spacelike, near any neighborhood of $0$.
The
proof in the general case, follows by
conjugating by
the exponential map $\exp_x$, where
$x$ is assumed to be, by contradiction,
a singular point of $V$.

Let's now show that the orbits of $V$ are geodesic.
As a Killing field, $V$ satisfies the following
anti-symmetry:
$<\nabla_VV,U>+<\nabla_UV,V>=0$,
for any vector field $U$.  Hence,  $<\nabla_VV, U> =
-(1/2)U.<V,V> = 0$, since $V$ is lightlike.
Therefore $\nabla_VV = 0$.

\end{proof}

The following proposition treats
lightlike Killing fields on
 constant curvature Lorentz manifolds.
Its  proof may be handled by a
standard  calculation.

\begin{proposition}
\label{light.field}
 Let ${\cal G} = {\cal G}(\tilde{N})$ be the Killing algebra of
$\tilde{N}$ (as above). Denote by
${\cal I} = {\cal I}(\tilde{N})$ the subset of $\cal G$ consisting
of lightlike Killing fields. Then:

i)  If $\tilde{N}$ is the anti de Sitter space (i.e. it
has negative curvature), then  ${\cal I}$ generates
${\cal G}$ as a vector-space.

ii) If $\tilde{N}$ is the Minkowski space (i.e. it
is flat), then the vector space generated by
${\cal I}$ equals the radical ${\bf R}^d$
of $\cal G$. More precisely, $X
\in {\cal I}$ if and only if $X$ is  parallel
(i.e. it generates a flow of translations) and is (somewhere)
isotropic.

iii) If $\tilde{N}$ is the de Sitter space, then
${\cal I} = \{0\}$.

\end{proposition}

\paragraph{Beginning.}
It is crucial to notice that,  if $X \in {\cal G}$ is lightlike,
as a Killing field on $\tilde{N}$) then
the same is true for $\bar{X}$,
as a Killing field on $\tilde{M}$.

Suppose by contradiction that $U \neq \tilde{M}$
(see Fact \ref{lorentzian} for notation), and let $N_1$
be the orbit of a point in the boundary of $U$.
By definition of $U$,
$N_1$ is lightlike (i.e.  is the metric on $TN_1$
is positive nondefinite).
If dim$N_1  \neq 0$, denote by
${\cal F}$ its  {\bf characteristic}
 foliation (of dimension 1),
i.e. that determined by the direction field
 $x \to N_1 \to (T_xN_1)^\perp \cap T_xN$.
 We denote by $Q$ the (local) quotient space $N_1/{\cal F}$.

%\paragraph{Lightlike submanifolds.}

\begin{fact}
If $X$ is a
lightlike Killing field, then the restriction
of $\bar{X}$ to $N_1$ is tangent to ${\cal F}$
(equivalently, the flow of such a Killing
field preserves individually the leaves of ${\cal F}$).
 Then  (from Proposition \ref{non.singular}),
we have:
 dim$N_1 \geq 1$, and the
 leaves of
$\cal F$ are lightlike geodesics (in $\tilde{M}$).
\end{fact}

\begin{proof}
A
lightlike field $X$
is tangent to ${\cal F}$, because
the direction of $\cal F$
is the unique isotropic direction tangent to
$N_1$.
\end{proof}

\paragraph{The anti de Sitter case.} In the case where  $\tilde{N}$ is the
 anti de Sitter space, $\cal I$ generates
$\cal G$, and therefore, $\cal G$ preserves (individually)
the leaves of
$\cal F$. Hence $N_1$ has dimension 1 (since
$N_1$ is a {\cal G}-orbit).
One then verifies that
$\cal G$
cannot act faithfully on such a manifold.
One  may see this in an easier way
in our situation here because $N_1$
reduces (at least
locally)
to the orbit of any element
$X \in {\cal G}$. Thus, it is an isotropic geodesic of $\tilde{M}$.
The action of $\cal G$ preserves the affine
parameter of this geodesic. Therefore,
$\cal G$ would  admit an injective homomorphism
in the affine group of
${\bf R}$, which is impossible.

\paragraph{The flat case.}  If
$N_1$ has dimension 1, we get a contradiction
as in the anti de Sitter case. If not,
consider the quotient space $Q = N_1 /{\cal F}$.
The ${\cal G}$-action on $N_1$,
factors through a faithful  action of $o(1, d-1)$
($= {\cal G}/{\bf R}^d$) on
$Q$.
Observe that
 $Q$  inherits a natural  Riemannian
metric. Indeed, the Lorentz metric restricted
to $N_1$ is positive degenerate, with kernel
$T{\cal F}$. But ${\cal F}$ is parameterized by
any lightlike Killing field $X \in {\cal I}$
(this is the meaning of the fact that the flow of  $\bar{X}$ preserves
individually
the leaves of ${\cal F}$). In particular,
the transversal action of the holonomy of ${\cal F}$
is equivalent to that of the flow of
$\bar{X}$. Therefore, it preserves the transverse metric,
which thus determines a metric on $Q$. This metric is invariant by
the $o(1, d-1)$-action. As in the
proof of Fact \ref{lorentzian},  since dim$Q \leq d-1$, we
have   dim$Q= d-1$, and furthermore,  $Q$ has constant curvature. Also, we
 recognize from
the list of Killing algebras of constant curvature manifolds that
$Q$ has constant negative curvature, i.e. $Q$
is a hyperbolic space.

It then follows that dim$N_1 = d$, and in particular that
%the foliation by
the orbits of $\cal G$ determine  a  regular
foliation near $N_1$. The idea,
to find a contradiction, is that, the leaves in $U_0$
have (intrinsic) constant 0 curvature, but not $N_1$,
because the quotient space $Q$ is hyperbolic.

To do this in
a more  rigorous manner, let $X_0 \in {\cal I}$ be a lightlike
Killing field,  and for all $x \in \tilde{M}$, consider
$Q_x $ the orbit  space of $X_0$ restricted to ${\cal G}x$.
% As above, in the case of $Q = Q_{x_1}$, $Q_x$ has a
%quotient Riemannian metric.
In a natural sense, the so obtained
quotient Riemannian spaces depend smoothly on $x \in \tilde{M}$.
 However, for $x \in U_0$, $Q_x$ is Euclidean (because
${\cal G}x$ is Minkowskian), but,
as stated above for $x \in N_1$, $Q_x$
is hyperbolic. This is a contradiction. Therefore $U =
\tilde{M}$.

\subsubsection{The de Sitter case.}
 In this case, ${\cal G} = o(1, d)$. Recall the example of
 the usual action
of $o(1, d)$ on the Minkowski space ${\bf R}^{1, d}$.
This
contrasts with
the case of nonpositive curvature,
 because the
orbits don't determine a (regular) foliation, since $0$ is a fixed point.
 Observe that there are
 3 types of  regular orbits: (Lorentzian) de Sitter orbits,
(Riemannian) hyperbolic orbits, and the isotropic cone at $0$
(without $0$),
which is lightlike.

In the general case, we have a partition of $\tilde{M}$ into
 degenerate and nondegenerate orbits.
 As it was  mentioned
in the  proof of Fact \ref{lorentzian}, a nondegenerate orbit
is a pseudo-Riemannian manifold
of dimension $\leq d$, and endowed with  a faithful
action
of $o(1, d)$, and hence, it is (locally and  up
to a multiple constant)
 either the (Riemannian) hyperbolic space or
the (Lorentz) de Sitter space.

 Therefore, we get a partition $D \cup H$ of the set of nondegenerate
orbits into, de Sitter and hyperbolic orbits, respectively. We write the
complementary set as  $\tilde{M}-D\cup H = L \cup S$, where $S$
is the subset of singular orbits, and
$L$ is the subset of degenerate but nonsingular orbits.

From Fact \ref{lorentzian}, for all
$x \in D$, $C_x$
contains the
isotropic cone of $T_x {\cal G}x$.
We have in fact equality:
$C_x =$
Cone$(T_x
{\cal G}x) $, by the
property
%choose of
of $d$ as a maximal dimension.

\paragraph{Invariance of the $\cal G$-foliation under the fundamental
group.}
 Now, our goal is   to show that essentially,
  the ${\cal G}$-foliation passes to
$M$. Indeed, consider $f \in \pi_1(M)$.
Denote by $  {\cal G}^\prime $ the image
by $f$
of the ${\cal G}$-action. It induces an analogous partition
$\tilde{M} = D^\prime \cup H^\prime \cup L^\prime \cup S^\prime$.
Of course, $C_{fx} = f(C_x)$, and hence
$C_y =$  Cone$(T_y  {\cal G}^\prime x)$,
for $y \in D^\prime = f(D)$.
In particular, in the open set
$D \cap D^\prime$, the $ {\cal G}$
and $  {\cal G}^\prime$ foliations coincide. Therefore,
by analyticity,
if $D \cap D^\prime \neq \emptyset$, the  $ {\cal G}$
and $  {\cal G}^\prime$-foliations coincide everywhere
in $\tilde{M}$.
Assume now that $D \cap D^\prime$ is empty and let
  $x \in D\cap H^\prime$. Consider the stabilizer ${\cal K}$
of $x$ for the $ {\cal G}^\prime$-action.
It is     isomorphic to $o(d)$
(the leaf $  {\cal G}^\prime x$ is isometric to
the hyperbolic space ${\bf H}^d$).
It
preserves $C_x$ and hence $T_x {\cal G}x$. Therefore,
we get a representation of $o(d)$ in $o(1, d-1)$
(the algebra of orthogonal transformations of
the Lorentz space $T_x  {\cal G}x$). But such a representation
must be trivial. Hence ${\cal K}$ acts trivially on
$T_x  {\cal G}x$, and thus also on its projection
on $T_x  {\cal G}^\prime x$. This projection is therefore
trivial, since the  ${\cal  K}$-action on
 $T_x  {\cal G}^\prime x$ is irreducible. This means that
the $  {\cal G}$ and $ {\cal G}^\prime$-orbits of $x$ are orthogonal at $x$
(if $x \in D \cap H^\prime$).
Thisextends by  analyticity to  all $\tilde{M}$.

In conclusion, exactly one of two  possibilities occurs
 for the  $  {\cal G}$ and $ {\cal G}^\prime$ foliations: they
coincide everywhere, or they are everywhere orthogonal.
But there is a finite number of mutually
orthogonal subspaces of dimension $d$ in a tangent space
$T_x\tilde{M}$. Therefore, there is a finite index subgroup
of
$\pi_1(M)$ which preserves the $ {\cal G}$-foliation.
For the sake of simplicity, we shall suppose that
$\pi_1(M)$ itself preserves the foliation.

\paragraph{Structure of orbits. The central flow.}
A    description of  orbits, similar to
that of the special case of the action
of $o(1, d)$ on the Minkowski space
${\bf R}^{1, d}$,
 holds in the
general case. Indeed, as  was  mentioned
above,
nondegenerate orbits are locally
isometric to
the de Sitter or to the hyperbolic space
of dimension $d$.
It remains to consider the case of lightlike
(nontrivial) orbits.
Let
$N_1$ be a such  orbit.  The quotient space $Q$ (of
the
characteristic
foliation of $N_1$, see above) is
a $o(1, d)$-homogeneous space of dimension $\leq d-1$.
A standard analysis
of subalgebras of
$o(1, d)$   shows
that $d-1$ is the minimal nontrivial dimension
of a space on which $o(1, d)$ acts. Furthermore,
the minimal dimension
is achieved
in the case of the usual
conformal action on the sphere $S^{d-1}$.

 Therefore, we have either  $Q $ is a single point, or
$Q $ is (locally) the conformal  sphere $S^{d-1}$.
  Of course, $Q$ cannot be  a single  point,
since otherwise, dim $N_1=1$,
but  $o(1, d)$ cannot act
nontrivially  on
${\bf R}$.
%would be reduced to
%a lightlike geodesic, endowed with its  affine action

Knowing that $Q$ is the usual conformal sphere, it is not difficult to
identify $N_1$ itself (locally) with
the isotropic cone of ${\bf R}^{1,d}$, as a $o(1, d)$-
homogeneous {\bf degenerate Riemannian space}.

%\paragraph{A contracting vector field tangent to the normal foliation.}
In the  Minkowski space ${\bf R}^{1, d}$, the action of
$o(1, d)$ on the isotropic cone commutes with the multiplication flow
$(t, x) \to e^t x$.
%In fact, up to a multiplicative constant, this
%is the unique flow which commutes with the $o(1, d)$-action.
In a similar way, one constructs a (local) {\bf central flow} $\tilde{\phi}^t$
on $D$
% (the subset of points with a degenerate, nonsingular leaf)
 which is tangent the $ {\cal G}$-foliation,
and commutes with the $  {\cal G}$-action. This flow passes
to a flow $\phi^t$ on the projection of $D$ in  $M$.

\paragraph{Structure near the singular set.} Let $x_0$ be a singular
point of the $ {\cal G}$-action. Then, we get an infinitesimal
representation of $o(1, d)$ in $o(1, n-1)$, the orthogonal algebra
of $T_{x_0}\tilde{M}$ ($n$ is the dimension of $M$).
 In a standard way, one may prove that such a representation
is equivalent to the usual inclusion $o(1, d) \subset o(1, n-1)$.
In particular,  there is  an orthogonal
decomposition $T_{x_0}M = E \bigoplus {\bf R}^{1, d}$, and
the infinitesimal action of
$o(1, d)$ is the product of the trivial action on $E$ and
the usual action of
$o(1, d)$ on ${\bf R}^{1, d}$.
The exponential map
$\exp_{x_0}$  conjugates (locally) the
action of $o(1, d)$ on $\tilde{M}$ with
its infinitesimal action on $T_{x_0}\tilde{M}$. In
particular the set of $ {\cal G}$-singular
points (near $x_0$) equals the geodesic spacelike
submanifold
$F= \exp_{x_0}E$. For $x \in F$,
the isotropic cone of $T_xF^\perp$ has
two sheets $Sh_x^\pm$ (here we don't mind
on the possibility of a continuous orientation
of these sheets when $x $ runs over
the singular set $S$).
%and $Sh_x^2$.
The degenerate
leaves
near $x_0$
are given by $\exp_xSh_x^\pm$
% or $\exp_xSh_x^2$
for $x \in F$. The (oriented non-parameterized) orbits of
$\tilde{\phi}^t$ (near $x_0$)
have the form $\exp_x tu$, $t >0$, where $u \in Sh_x^\pm$,
and $x \in F$.
Here, it is essential to observe that
$F$ is a repulsor of $\tilde{\phi}^t$. More precisely,
for
$y \in L$, near $x_0$, the orbit
$\tilde{\phi}^t(y)$ converges to a point $x \in F$,
when
$t \to -\infty$.

From this, one sees that the flow
$\tilde{\phi}^t$
can be continuously extended to the $ {\cal G}$-singular set $S$,
by letting it act trivially on
$S$. Therefore, the flow $\tilde{\phi}^t$
is now defined on the  closed subset
$D \cup S$. Its quotient flow
$\phi^t$ is thus defined on a compact space and in
particular
$\tilde{\phi}^t$ is complete.

Let $V$ be a ``conical'' neighborhood of
$S$ in $S \cup L$, that is,
$V$ has the form
$V = \cup_{x \in S} V_x$,
where
$V_x$ is the intersection
of the isotropic cone
of $T_xS^\perp$, with
a ball
in $T_x\tilde{S}$ (with
respect to any Riemannian metric
on $\tilde{M}$). Suppose that
$V$ is $\pi_1$-invariant (choose
$V$ to come from
a similar neighborhood of the projection of
$S$
in $M$). From the repelling property described above,
the complement $L^*= L-V$ is $\tilde{\phi}^t$-invariant
for
$t>0$.

%\paragraph{Contradiction.}
% in the compact manifold.}
 We will get a contradiction by showing that
$D$ is (a codimension 1) submanifold of $\tilde{M}$,
and for $t$ big enough, the Jacobian
$Jac(D_x \tilde{\phi}^t)$ is uniformly
$>1$, for $x \in D^*$.

\paragraph{Structure of $L$.}  From
above, the singular set $S$ has  codimension at least
$3$, so
$\tilde{M} -S$ ($ = D \cup L \cup H$)
is connected. There, the orthogonal plane field
of the $ {\cal G}$-foliation
has constant dimension $n-d$. Therefore, it
is analytic, and integrable, and has geodesic leaves (since
this is the case
in an  open subset of $D$).
We denote by $\cal L$ its tangent foliation (defined on
$\tilde{M}-S$).

For  $x \in L$,  we have $T_x{\cal L} =
 (T_x {\cal G}x)^\perp$. Hence, as
the $ {\cal G}$-orbit of $x$, the
geodesic
leaf
${\cal L}_x$ is degenerate. Since it is geodesic,
$T_y {\cal L}_x$
is degenerate, for all $y\in {\cal L}_x$, and thus
the $ {\cal G}$-orbit of
$y$ is degenerate, that is $y \in L$. In conclusion,
$L$
is invariant by both
 the $ {\cal G}$
and the
$\cal L$-foliations.

For $x \in L$, the intersection $T_x{\cal L} \cap
 T_x {\cal G}x$, is exactly the
isotropic direction
of $
 T_x {\cal G}_x$, which is
nothing but the tangent direction
of $\tilde{\phi}^t$
at $x$. In other words, the leaves
$ {\cal G}x$ and ${\cal L}_x$
meet along
the $\tilde{\phi}^t$-orbit of $x$.

From this, we infer that, around any $x $ in $ \tilde{M}$,
there
is a hypersurface (of $\tilde{M}$) contained in $L$. In
fact, $L$ which is an analytic set of
$\tilde{M}-S$ is a regular hypersurface.  To see this, one
constructs  an analytic distribution $\Delta$, for
which $L$ is a union of integral leaves, as follows.
Locally, $\Delta$ is generated by a family
 $X_1, \ldots, X_d, Y_1, \ldots, Y_{n-d}$
of vector fields, where
 $X_1, \ldots, X_d$ (resp.  $Y_1, \ldots, Y_{n-d}$)
 generate
the tangent space of the
$ {\cal G}$-foliation (resp. the tangent space of the
foliation
$\cal L$). From the previous discussion,
 one sees that  $L$
is the union of singular leaves of $\Delta$.

Let  $g$ be a $\pi_1$-invariant Riemannian
metric on
$\tilde{M}$. Denote by
$X$ the vector field generating the
flow $\tilde{\phi}$, and for $x \in L$, let
$E^u(x) $ (resp. $E^0(x)$) be the orthogonal
(with respect to
$g$)
of $X(x)$ in $T_x {\cal G}x$ (resp.
$T_x{\cal L}_x$). These two plane
fields on $L$ are spacelike (with respect
to
the Lorentz metric of $\tilde{M}$). We change
the Riemannian metric on
$L$,  by equipping $E^u$ and
$E^0$ with the restriction of Lorentz metric,
and  decreeing
that $X$, $E^u$ and $E^0$
are orthogonal.
 With respect to this
new Riemannian
metric, we have the following
relations: $ \vert D\tilde{\phi}^t u \vert = e^t\vert u\vert$,
for $u \in E^u$ and $ \vert D\tilde{\phi}^t u \vert = \vert u\vert$,
for
$u \in E^0$ (and of course $ \vert D\tilde{\phi}^t u \vert = \vert u\vert$,
for
$u \in {\bf R}X$).
The first relation comes from the fact  that
this is the case, in the model space (of
$ {\cal G}x$, which is the isotropic cone
of
the Minkowski space
${\bf R}^{1, d}$, see above).

The second relation follows from the following
general property of
lightlike {\bf geodesic}
submanifolds
in Lorentz manifold \cite{Zeg.geodesic.foliations}. It is that
 the (one dimensional)
 characteristic (isotropic)  foliation of such a lightlike
submanifold is a
transversely Riemannian
foliation.  A flow parameterizing
the characteristic  foliation, preserves the degenerate
Riemannian
metric of the lightlike geodesic submanifold.

 Now, the projection of
$D^*$ in $M$ is a compact  manifold (with boundary), preserved
by the (semi-)flow $\phi^t$ ($t >0$), and
$Jac(\phi^t) = e^{(d-1)t}$. This is impossible. This means that
the lightlike locus $L$ is empty, and therefore,
$\tilde{M} = D$, that is all the $  {\cal G}$-leaves
are of de Sitter type. Thus, as desired,   $M$
is  everywhere, locally a warped product.

%in a neighborhoud of point $x$ of $L$. $T_x{ {\cal G}x$ (rep. $T_

\subsection{Completeness}

\paragraph{Global  topological product.}
The  foliation $ {\cal G}^\perp$ is
  geodesic, in the sense of the Lorentz
metric. Therefore, the holonomy
of $ {\cal G}$ preserves the metric on
$T {\cal G}^\perp$.
Thus $ {\cal G}$
 is a transversely Riemannian foliation, since
$ {\cal G}^\perp$
is spacelike.

Transform the Lorentz metric of $M$ into a
 Riemannian metric, by keeping the same metric
on $T {\cal G}^\perp$, keeping $ {\cal G}$ and
$ {\cal G}^\perp$ orthogonal, and
choosing any ($\pi_1$-invariant)  Riemannian metric on
$T {\cal G}$.

The foliation $ {\cal G}$
 is still  transversely Riemannian in the sense of the new metric,
 because we have not changed the metric on its orthogonal.
It then follows that
${\cal G}^\perp$
is geodesic (in the sense of the new Riemannian metric).

Now, we use a result of \cite{Hebda}
which states that if  a geodesic foliation
${\cal L}$ of a compact Riemannian manifold $M$
admits an orthogonal foliation
$\cal  N$ (that is $T{\cal L}^\perp$ is integrable),
then the pair $(\tilde{\cal L}, \tilde{\cal N})$
gives a global topological product
in $\tilde{M}$.  More precisely,
let $\tilde{\cal L}_x$
and $\tilde{\cal N}_x$ be the leaves of a point
$x $ of $\tilde{M}$. Then, the
inclusion of $\tilde{\cal L}_x \cup \tilde{\cal N}_x$
 in $\tilde{M}$ extends to a diffeomorphism
$\tilde{\cal L}_x \times \tilde{\cal N}_x   \to \tilde{M}$
 sending the foliations determined
by the factors,  to
$\cal L$ and $\cal N$,
respectively.

\paragraph{A new Lorentz product metric.}
Let $x_0 \in \tilde{M}$, then   $\tilde{M} $
is homeomorphic to
$\tilde{\cal L}_{x_0} \times \tilde{\cal N}_{x_0}  $,
endowed with a warped
metric $h \bigoplus w g$, where
$h$ (resp. $g$) is the Riemannian (resp. Lorentzian) metric
on
$\tilde{\cal L}_{x_0} $
(resp.
$\tilde{\cal N}_{x_0}   $), and
$w:
 \tilde{\cal L}_{x_0}   \to {\bf R}$ is
a warping function.

Our aim here, is to show that
the product metric $h \bigoplus g$ is
$\pi_1(M)$-invariant (and hence descend to a
locally product metric on
$M$).  This is equivalent
to the invariance of
the warping function $w$ by the $\pi_1$-action
on
$\tilde{\cal L}_{x_0}$.

If $ \tilde{\cal N}_{x_0}   $
%In the case where the factor $N$
is not flat,
   $w$ is
$\pi_1$-invariant, since
$w(x) = \kappa(x)^{-2}$,
where $\kappa(x)$ is
the sectional curvature
  of $\tilde{\cal N}_x$.

Now, we  give briefly  an  idea of how to prove the invariance
of
$w$ in the flat case. Let's start
by considering the case
dim
$\tilde{\cal L}_{x_0}= 1$. In $M$, we have a codimension
one transversely Riemannian foliation
${\cal N}$. Let $\phi^t$ be the flow
generated by a unit vector field orthogonal
to
${\cal N}$. Then,  $\phi^t$  preserves ${\cal N}$ (this
 is exactly the meaning of
${\cal N}$
being transversely Riemannian).  In fact,
 because  of the local warped product
structure,  $\phi^t$
sends  leaves homotheticaly to  leaves.  More precisely,
let $\tilde{\phi}^t$
be the lift of $\phi^t$ to  $\tilde{M}$. Then,
for $x$, $y$ $\in \tilde{\cal L}_{x_0}$,
 we have
 Jac$( \tilde{\phi}^t_x) =$
 Jac$( \tilde{\phi}^t_x \vert \tilde{\cal N}) =
 (w(y)/w(x))^{n-1/2}$, where $t$ is such that,
 $ \tilde{\phi}^t(x) = y$ ($n =$ dim$M$).

There are two possibilities. The first is that
 all the leaves of ${\cal N}$
are dense. In this case, Jac$({\phi}^t) = $ constant, and hence equals 1,
and thus,  $w $  is constant.
The second  case is that  all the leaves of ${\cal N}$
are closed. We then have: for $x$, $y$ $\in
\tilde{\cal L}_{x_0}$,
Vol${\cal N}_{\pi(x)}/$
Vol${\cal N}_{\pi(y)} = (w(x)/w(y))^{n-1)/2}$, where
$\pi: \tilde {M} \to M$
is the projection. Therefore, $w$ is
$\pi_1$-invariant.

The same proof work in the higher
dimensional case if we suppose that there
is a parallelism, i.e. a frame of vector
fields $X^1, \ldots, X^k$, tangent to
${\cal L}$, preserving the foliation ${\cal N}$, and
also
preserving
the volume along
${\cal L}$.  But, such a parallelism
 is exactly  what Molino's theory on Riemannian
foliations \cite{Mol} yields, up to passing to
 another foliation naturally associated to
${\cal N}$.

\paragraph{Completeness along  the constant curvature factor.}
Recall that compact Lorentz
manifolds of constant   curvature curvature, are
(geodesically) complete.   This result
 was proved by  Y. Carri\`ere \cite{Car}, in the flat case,
 and then the proof was
adapted to the general case by
B. Klingler \cite{Kli}. Our observation  here is that the Carri\`ere's
proof  may be easily updated to handle the case
of compact Lorentz manifolds whose the universal
cover is a global  (direct)  product of a Lorentz
manifold of constant curvature by a Riemannian manifold.
The point
is to develop,  $\tilde{M}$ which is the   product
$\tilde{\cal L}_{x_0} \times \tilde{\cal N}_{x_0} $,
into the
product
$ \overline{\tilde{M}} =
\tilde{\cal L}_{x_0} \times \tilde{ N} $
  where
$\tilde{N}$ is the simply connected
Lorentz manifold with the same curvature as
$\tilde{\cal N}_{x_0}$.
  Then, instead of triangles, as used
in the Carri\`ere's proof, we use subsets
of
$\overline{\tilde{M}}$ of the form $B \times \Delta$,
where
$B$ is a ball of
$\tilde{\cal L}_{x_0}$
and
  $\Delta$
 is a triangle in
$\tilde{N}$.
This leads to the conclusion
that
$\tilde{M}$
is isometric
to
$\overline{\tilde{M}}$.

\paragraph{Completeness.} We infer from
\cite{Shanz}, that $M$, endowed with the old warped
Lorentz metric, is complete. Indeed
the warping function $w$ is  bounded, since it is
$\pi_1$-invariant.

\paragraph{End of the proof of Theorem \protect\ref{global}.}  It
 remains to
show
that the
% constant curvature factor
 factor $\tilde{N}$ has nonpositive (constant) curvature. This fact
was observed in Part I, \S 15, as being a slight generalization
of the Calabi-Markus
phenomena.

\section{Analytic bi-polarized manifolds. Proof of the second half of
 Theorem \protect\ref{classification}}

The possibilities  (1) and (2) of Theorem
\ref{classification}, correspond to
the cases ${\cal W}(M)
\neq \emptyset$, and
 ${\cal W}(M)
= \emptyset$, respectively. The structure of $M$ when
 ${\cal W}(M)
\neq \emptyset$ follows from Theorem \ref{global}. The
goal of the present section , is to study the other situation.
 So, let
$M$  be  a compact Lorentz manifold, with
Isom$M$ noncompact, and ${\cal W}(M) =
\emptyset$.  We will show that,
in the analytic case,   $M$
satisfies the description stated in point (2) of Theorem
\ref{classification}.
As  was said in the introduction, the noncompactness
 of Isom$M$ implies that
card $C_x \geq 1$ for
all $x \in M$ (this follows from
Part I, or a direct proof).

\paragraph{Bi-polarization.}
%In the sequel, except for the following  proposition,
In what follows, except for Proposition \ref{local.polarization},
$M$ will
be supposed to be analytic. In fact,  in the
proof of this proposition itself, we
will use the fact that $M_{smooth}$
 is dense in $M$.  This was observed to
be true in the analytical case (see  \ref{hyper.warped.analytic}),  but
its proof in the smooth case is harder (\cite{Zeg.tautologic}).
Also, from
 the  latter  reference, we infer that one may slightly change
the definition of
$M_{smooth}$
so that it remains open and dense, and
not only the map $x \to E_x$ (the space generated by
$C_x$) is smooth, but the map $x \to C_x$
itself is semi-continuous in an obvious manner. Again,
in the analytic case, this fact follows from that the integrability
 domain
is an analytic set.

\begin{proposition}
\label{local.polarization}
 Let $M$  be  a compact Lorentz manifold, with
Isom$M$ noncompact, and ${\cal W}(M) =
\emptyset$. Then,
$1 \leq$ card $ C_x \leq 2$,
for
all $x \in M_{smooth}$.

In particular $M$ is bi-polarized, that is
Isom$(M)$ preserves  two (perhaps  identical)
lightlike geodesic foliations ${\cal F}_1$ and
${\cal F}_2$
(see \cite{Part.I}, \S 11, for details).

In fact, in $M_{smooth}$, we have
$C_x =  \{ T_x{\cal F}_1^\perp, T_x{\cal F}_2^\perp\}$.

%In particular
\end{proposition}

\begin{proof}
 From  Theorem \ref{hyper.warped.2}, we infer,
since ${\cal W}(M) = \emptyset$,
that int $(M_{\leq \hbox{dim} M})$ is dense in $M_{smooth}$.
%We claim that: int($M_{\leq \hbox{dim} M}) \subset
%M_{\leq 2}$.
 Consider $A=  \cup \{C_x\; / x \in \hbox{int}(M_{\leq
\hbox{dim} M})\}$.
This is an Isom$M$-invariant  subset of the
projective isotropic tangent  bundle ${\bf P}T^0M$,
with  finite  fibers $A_x$ (over $M$). Furthermore, it is measurable,
since $x \to C_x$ is semi-continuous \cite{Zeg.tautologic}. From
the barycenter construction (Theorem 2.6 of Part I), Isom$M$
 would be compact,
if we didn't have card$A_x \leq 2$ almost everywhere. This
implies that
almost everywhere, dim$E_x \leq 2$. This extends to all $M_{smooth}$,
since the map  $x \to E_x$
is smooth, and hence  card $C_x \leq 2$, for all
$x \in M_{smooth}$.

In particular, the
fibers of the limit set $\Lambda$ of Isom$M$ in ${\bf P}T^0M$
%is a union of continuous (in fact Lipschitz) sections
%(which determine lightlike geodesic foliations). Thus, everywhere
%in $M$, the fiber
also satisfy : card $\Lambda_x \leq 2$.
 %contains at most 2 elements.
Therefore
  Isom$M$ is elementary, or equivalently, $M$ is bi-polarized
(see Part I, \S 11).

Finally, let's check the
equality $C_x =  \{ T_x{\cal F}_1^\perp, T_x{\cal F}_2^\perp\}$,
along $M_{smooth}$.
%, follows from the continuity
We have seen that,  dim$E_x \leq 2$, for $x \in M_{smooth}$.
 Therefore, if $x $ belongs to the transversality set
 $ T =
  \{ x \in M / T_x{\cal F}_1^\perp
 \neq T_x{\cal F}_2^\perp\}$, then,
 card$ \{ T_x{\cal F}_1^\perp, T_x{\cal F}_2^\perp\}= 2$, and
hence $C_x =  \{ T_x{\cal F}_1^\perp, T_x{\cal F}_2^\perp\}$.
This extends by continuity to the closure of $T$
in $M_{smooth}$.

In the coincidence set $C= M-T$, we have:
card$\{ T_x{\cal F}_1^\perp, T_x{\cal F}_2^\perp\}= 1$.
Suppose that in some (open) component $U$  of int($M_{smooth} \cap C$),
we have card$C_x = 2$. Then, we get a continuous  Isom$M$-invariant
 isotropic direction field
$X$ along $U$, such
that
$C_x= \{ T_x{\cal F}_1^\perp = T_x{\cal F}_2^\perp , X(x)\}$.
 This contradicts the ``north-south''
dynamical behavior of the derivative action of
Isom$M$ on ${\bf P}T^0M$. Indeed, on $U$, the
dynamics has only one ``pole'' $ T{\cal F}_1^\perp \vert U
 = T{\cal F}_2^\perp \vert U$.
Thus, for an Isom$M$-recurrent point $x$ (which  exists since $U$
is open), we must have
$X(x)= T_x{\cal F}_1^\perp  = T_x{\cal F}_2^\perp$, which contradicts
the definition of $X$.
This finishes the proof of the proposition.
\end{proof}

\paragraph{Local bi-polarization in the analytic case.}
The foliations ${\cal F}_1$ and ${\cal F}_2$
determine two sections of the Grassmann
bundle $Gr^0 (M) \to M$.  Their
images $P({\cal F}_1)$ and $P({\cal F}_2)$
are
 (topological) submanifolds of $Gr^0(M)$,  homeomorphic to $M$.
Let  ${\cal D}^0$ denote $   P({\cal F}_1) \cup P({\cal F}_2)$.

Consider ${\cal D}$,  the integrability domain of $\tau$ on
$Gr^0(M)$.
Obviously,  ${\cal D}^0  \subset {\cal D}$.

\begin{fact}  Keep the hypotheses of the
 previous proposition, and assume
that
$M$ is analytic. Let
$p \in {\cal D}- {\cal D}^0$,
and $ \Sigma \subset Gr^0(M)$, a small
tranversal to the tautological plane field
$\tau$, containing $p$.
Then $p$ is isolated in $\Sigma \cap {\cal D}$.
\end{fact}

\begin{proof}
Suppose the contrary. Then, by
analyticity, $ \Sigma \cap {\cal D}$ contains a curve,
which is in fact contained in ${\cal D}-{\cal D}^0$ if
$\Sigma$ is small enough (since ${\cal D}-{\cal D}^0$
is open in ${\cal D}$). In other words, we obtain
a one parameter family of lightlike geodesic
hypersurfaces, which are not leaves of ${\cal F}_1$
or ${\cal F}_2$.  This one parameter family fills,
at least, some open set $U$, say.
This contradicts the validity
of the equality $C_x =  \{ T_x{\cal F}_1^\perp, T_x{\cal F}_2^\perp\}$,
in the open dense set $M_{smooth}$.
\end{proof}

\begin{corollary}
%[Local bi-polarization]
%Keep the hypotheses of the previous proposition, and assume
%that
%$M$ is analytic. Then:
 We have: ${\cal D} = {\cal D}^0$ (Equivalently:
$C_x = \{ T_x{\cal F}_1^\perp, T_x{\cal F}_2^\perp\}$,
 for all $x \in M$).

\end{corollary}

\begin{proof}
We deduce from the fact above that ${\cal D}-{\cal D}^0$
is closed in ${\cal D}$. Indeed, let $p \in {\cal D}^0$,  and
$\Sigma$ a transversal as  above. Then,
from the previous fact, $\Sigma \cap ({\cal D}-{\cal D}^0)$
is open and discrete in the analytic set
$\Sigma \cap {\cal D}$. This latter set has finitely
many connected components. Therefore,
$\Sigma \cap ({\cal D}-{\cal D}^0)$ must be finite. In particular
$p$ cannot be an accumulation point of ${\cal D}-{\cal D}^0$,
so  $ {\cal D}-{\cal D}^0$ is closed.

It then follows that ${\cal D}-{\cal D}^0$
consists of a finite union of closed leaves of $\tau$.
By projecting  in $M$, we get
 closed lightlike geodesic hypersurfaces
$S^1, \ldots S^k$.
Up to a subgroup of finite index, we may suppose
that Isom$(M)$ preserves each hypersurface $S^i$.
 Observe that $S^i$ is nowhere tangent
to ${\cal F}_1$ (or ${\cal F}_2)$)
since otherwise, it would be a leaf of this foliation.

As in the proof of the proposition above,
for $f \in \hbox{Isom}(M)$,
the action of $Df$ on the projective
isotropic cone ${\bf P}T^0M \vert S^i$,
 along $S^i$,
preserves $(TS^i)^\perp$, which contradicts
 the fact that it is  a north-south
%property of this
dynamics determined by
the attractor-repulsor pair $((T{\cal F}_1)^\perp
\vert S^i, (T{\cal F}_2)^\perp
\vert S ^i)$.

\end{proof}

\paragraph{Regularity.} The analytic set ${\cal D}$
equals $P({\cal F}_1) \cup P({\cal F}_2)$,
 a union of two topological manifolds.
As above, let  $\Sigma$ be a transversal to $\tau$, at a point
$p \in {\cal D}$, and consider the 1-dimensional
(local) analytic set $A= {\cal D} \cap \Sigma$. Topologically,
$A$ is a union of one or two topological curves, depending on
 the projection of $p$ in  $M$
belongs to the transversality set $T$, or
to the coincidence set $C$ respectively (recall that
$C$ is the closed subset of $M$ where
${\cal F}_1$ and ${\cal F}_2$
are tangent).

From the structure theory of
1-dimensional (real) analytic sets (\cite{Chir}),
$A$ is a union finitely many branches, i.e.
 images of  analytic curves. Here, an
analytic curve means an analytic injective map
from an interval $]-\epsilon, + \epsilon [$
 to $Gr^0(M)$, and sending 0 to
$p$. In our situation, there are  two analytic
curves $c_1$ and $c_2$, such that
$A = \hbox{Image}(c_1)$ if $p$ projects on $T$, and
$A = \hbox{Image}(c_1) \cup \hbox{Image}(c_2)$
 if $p$ projects on $C$.

One starts (if necessary) by permuting branches,
or equivalently changing the
foliations ${\cal F}_1$ and ${\cal F}_2$, so
that, above  the coincidence set $C$, we
have: $\Sigma \cap P({\cal F}_1) = \hbox{Image}(c_1)$
and $\Sigma \cap P({\cal F}_2) = \hbox{Image}(c_2)$.
One sees that this manipulation gives rise to new foliations
(always denoted
${\cal F}_1$ and
${\cal F}_2$) which coincide with the former ones in each side
of the coincidence set $C$.  Observe  furthermore that
the new decomposition
${\cal D}= P({\cal F}_1) \cup P({\cal F}_2)$ is
canonical, since  the decomposition
$A = \hbox{Image}(c_1) \cup \hbox{Image}(c_2)$
is canonical
(although
  there
are no canonical  parameterizations
$c_1$ and $c_2$).

Now, let's consider one of the foliations, say ${\cal F}_1$,
and perform some change of the induced
analytic structure of $P({\cal F}_1)$, so
that it becomes an analytic manifold.  Consider
as above a transversal set $A = \Sigma \cap P({\cal F}_1)
= \hbox{Image}(c_1)$. Then, $p$ is a regular point
of $P({\cal F}_1)$, if and only if $p$
 is a regular point of $A$, if and only if we may
choose $c_1$
 having 0 (the pre-image of $p$)
as an immersion point.

If $p$ is a singular point, we choose a less singular curve
$c_1$, and define an (abstract) analytic structure
on $A$ such that
$c_1$ is a parameterization (i.e.
$c_1^{-1}$ is a  chart).  One then changes  the analytic
structure of $P({\cal F}_1)$
near $p$ accordingly.
 More precisely,  in $P({\cal F}_1)$, $p$
 admits a neighborhood which is the image
of a foliated homeomorphic analytic map
$\phi: (t, u) \in I \times {\bf R}^{n-1} \to
(c_1(t), f(u))  \in  P({\cal F}_1)$,
where $f$ is an analytic diffeomorphism. Then
we decree
that $\phi^{-1}$
 is a chart for the new analytic structure. One
then can verify that this is
a well defined analytic structure on
$P({\cal F}_1)$,  for which
${\cal F}_1$
(or equivalently the plane field
$\tau$)
 is analytic. Endow $M$
with
the pull-back of the so defined  structure
on $P({\cal F}_1)$
by the section map $x \to T_x{\cal F}_1$. Then
${\cal F}_1$ is analytic with respect to this structure.
This finishes the proof
of the second  half of Theorem \ref{classification}.

%\end{proof}

UMPA ENS Lyon, 46, all\'ee d'Italie \\
69364 Lyon cedex 7 FRANCE\\
Zeghib@umpa.ens-lyon.fr

\end{document}